# VARIABLE SELECTION USING MM ALGORITHMS

By David R. Hunter and Runze Li[1]

*Pennsylvania State University*

Variable selection is fundamental to high-dimensional statistical modeling. Many variable selection techniques may be implemented by maximum penalized likelihood using various penalty functions. Optimizing the penalized likelihood function is often challenging because it may be nondifferentiable and/or nonconcave. This article proposes a new class of algorithms for finding a maximizer of the penalized likelihood for a broad class of penalty functions. These algorithms operate by perturbing the penalty function slightly to render it differentiable, then optimizing this differentiable function using a minorize–maximize (MM) algorithm. MM algorithms are useful extensions of the well-known class of EM algorithms, a fact that allows us to analyze the local and global convergence of the proposed algorithm using some of the techniques employed for EM algorithms. In particular, we prove that when our MM algorithms converge, they must converge to a desirable point; we also discuss conditions under which this convergence may be guaranteed. We exploit the Newton–Raphson-like aspect of these algorithms to propose a sandwich estimator for the standard errors of the estimators. Our method performs well in numerical tests.

**1. Introduction.** Fan and Li [7] discuss a family of variable selection methods that adopt a penalized likelihood approach. This family includes well-established methods such as AIC and BIC, as well as more recent methods such as bridge regression [11], LASSO [23] and SCAD [2, 7]. What all of these methods share is the fact that they require the maximization of a penalized likelihood function. Even when the log-likelihood itself is relatively easy to maximize, the penalized version may present numerical challenges. For example, in the case of SCAD or LASSO or bridge regression, the penalized log-likelihood function is nondifferentiable; with SCAD or bridge

Received May 2003; revised November 2004.
[1]Supported by NSF Grants DMS-03-48869 and CCF-04-30349, and NIH Grant NIDA 2-P50-DA-10075.
*AMS 2000 subject classifications.* 62J12, 65C20.
*Key words and phrases.* AIC, BIC, EM algorithm, LASSO, MM algorithm, penalized likelihood, oracle estimator, SCAD.

regression, the function is also nonconcave. To perform the maximization, Fan and Li [7] propose a new and generic algorithm based on local quadratic approximation (LQA). In this article we demonstrate and explore a connection between the LQA algorithm and minorization–maximization (MM) algorithms [14], which shows that many different variable selection techniques may be accomplished using the same algorithmic techniques.

MM algorithms exploit an optimization technique that extends the central idea of EM algorithms [6] to situations not necessarily involving missing data nor even maximum likelihood estimation. The connection between LQA and MM enables us to analyze the convergence of the local quadratic approximation algorithm using techniques related to EM algorithms [16, 19, 20, 24]. Furthermore, we extend the local quadratic approximation idea here by forming a slightly perturbed objective function to maximize. This perturbation solves two problems at once. First, it renders the objective function differentiable, which allows us to prove results regarding the convergence of the MM algorithms discussed here. Second, it repairs one of the drawbacks that the LQA algorithm shares with forward variable selection: namely, if a covariate is deleted at any step in the LQA algorithm, it will necessarily be excluded from the final selected model. We discuss how to decide a priori how large a perturbation to choose when implementing this method and make specific comments about the price paid for using this perturbation.

The new algorithm we propose retains virtues of the Newton–Raphson algorithm, which among other things allows us to compute a standard error for the resulting estimator via a sandwich formula. It is also numerically stable and is never forced to delete a covariate permanently in the process of iteration. The general convergence results known for MM algorithms imply among other things that the newly proposed algorithm converges correctly to the maximizer of the perturbed penalized likelihood whenever this maximizer is the unique local maximum. The linear rate of convergence of the algorithm is governed by the largest eigenvalue of the derivative of the algorithm map.

The rest of the article is organized as follows. Section 2 briefly introduces the variable selection problem and the penalized likelihood approach. After providing some background on MM algorithms, Section 3 explains their connection to the LQA idea, then provides a modification to LQA that may be shown to be an MM algorithm for maximizing a perturbed version of the penalized likelihood. Various convergence properties of this new MM algorithm are also covered in Section 3. Section 4 describes a method of estimating covariances and presents numerical tests of the algorithm on a set of four diverse problems. Finally, Section 5 discusses the numerical results and offers some broad comparisons among the competing methods studied in Section 4. Some proofs appear in the Appendix.

**2. Variable selection via maximum penalized likelihood.** Suppose that $\{(\mathbf{x}_i, y_i): i = 1, \ldots, n\}$ is a random sample with conditional log-likelihood $\ell_i(\boldsymbol{\beta}, \boldsymbol{\phi}) \equiv \ell_i(\mathbf{x}_i^T\boldsymbol{\beta}, y_i, \boldsymbol{\phi})$ given $\mathbf{x}_i$. Typically, the $y_i$ are response variables that depend on the predictors $\mathbf{x}_i$ through a linear combination $\mathbf{x}_i^T\boldsymbol{\beta}$, and $\boldsymbol{\phi}$ is a dispersion parameter. Some of the components of $\boldsymbol{\beta}$ may be zero, which means that the corresponding predictors do not influence the response. The goal of variable selection is to identify those components of $\boldsymbol{\beta}$ that are zero. A secondary goal in this article will be to estimate the nonzero components of $\boldsymbol{\beta}$.

In some variable selection applications, such as standard Poisson or logistic regression, no dispersion parameter $\boldsymbol{\phi}$ exists. In other applications, such as linear regression, $\boldsymbol{\phi}$ is to be estimated separately after $\boldsymbol{\beta}$ is estimated. Therefore, the penalized likelihood approach does not penalize $\boldsymbol{\phi}$, so we simplify notation in the remainder of this article by eliminating explicit reference to $\boldsymbol{\phi}$. In particular, $\ell_i(\boldsymbol{\beta}, \boldsymbol{\phi})$ will be written $\ell_i(\boldsymbol{\beta})$. This is standard practice in the variable selection literature; see, for example, [7, 11, 21, 23].

Many variable selection criteria arise as special cases of the general formulation discussed in [7], where the penalized likelihood function takes the form

$$Q(\boldsymbol{\beta}) = \sum_{i=1}^{n} \ell_i(\boldsymbol{\beta}) - n\sum_{j=1}^{d} \lambda_j p_j(|\beta_j|) \equiv \ell(\boldsymbol{\beta}) - n\sum_{j=1}^{d} \lambda_j p_j(|\beta_j|). \tag{2.1}$$

In (2.1), the $p_j(\cdot)$ are given nonnegative penalty functions, $d$ is the dimension of the covariate vector $\mathbf{x}_i$ and the $\lambda_j$ are tuning parameters controlling model complexity. The selected model based on the maximized penalized likelihood (2.1) satisfies $\beta_j = 0$ for certain $\beta_j$'s, which accordingly are not included in this final model, and so model estimation is performed at the same time as model selection. Often the $\lambda_j$ may be chosen by a data-driven approach such as cross-validation or generalized cross-validation [5].

The penalty functions $p_j(\cdot)$ and the tuning parameters $\lambda_j$ are not necessarily the same for all $j$. This allows one to incorporate hierarchical prior information for the unknown coefficients by using different penalty functions and taking different values of $\lambda_j$ for the different regression coefficients. For instance, one may not be willing to penalize important factors in practice. For ease of presentation, we assume throughout this article that the same penalization is applied to every component of $\boldsymbol{\beta}$ and write $\lambda_j p_j(|\beta_j|)$ as $p_\lambda(|\beta_j|)$, which implies that the penalty function is allowed to depend on $\lambda$. Extensions to situations with different penalty functions for each component of $\boldsymbol{\beta}$ do not involve any extra difficulties except more tedious notation.

Many well-known variable selection criteria are special cases of the penalized likelihood of (2.1). For instance, consider the $L_0$ penalty $p_\lambda(|\beta|) = 0.5\lambda^2 \times I(|\beta| \neq 0)$, also called the entropy penalty in the literature, where $I(\cdot)$ is an indicator function. Note that the dimension or the size of a model equals $\sum_j I(|\beta_j| \neq 0)$, the number of nonzero regression coefficients in the model. In other words, the penalized likelihood (2.1) with the entropy penalty can be rewritten as

$$\ell(\boldsymbol{\beta}) - 0.5n\lambda^2|M|,$$

where $|M| = \sum_j I(|\beta_j| \neq 0)$, the size of the underlying candidate model. Hence, several popular variable selection criteria can be derived from the penalized likelihood (2.1) by choosing different values of $\lambda$. For instance, the AIC (or $C_p$) and BIC criteria correspond to $\lambda = \sqrt{2/n}$ and $\sqrt{(\log n)/n}$, respectively, although these criteria were motivated from different principles. Similar in its effect to the entropy penalty function is the hard thresholding penalty function (see [1]) given by

$$p_\lambda(|\beta|) = \lambda^2 - (|\beta| - \lambda)^2 I(|\beta| < \lambda).$$

Recently, many authors have been working on penalized least squares with the $L_q$ penalty $p_\lambda(|\beta|) = \lambda|\beta|^q$. Indeed, bridge regression is the solution of penalized least squares with the $L_q$ penalty [11]. It is well known that ridge regression is the solution of penalized likelihood with the $L_2$ penalty. The $L_1$ penalty results in LASSO, proposed by Tibshirani [23]. Finally, there is the smoothly clipped absolute deviation (SCAD) penalty of Fan and Li [7]. For fixed $a > 2$, the SCAD penalty is the continuous function $p_\lambda(\cdot)$ defined by $p_\lambda(0) = 0$ and, for $\beta \neq 0$,

$$p'_\lambda(|\beta|) = \lambda I(|\beta| \leq \lambda) + \frac{(a\lambda - |\beta|)_+}{a-1} I(|\beta| > \lambda), \tag{2.2}$$

where throughout this article $p'_\lambda(|\beta|)$ denotes the derivative of $p_\lambda(\cdot)$ evaluated at $|\beta|$.

Letting $p'_\lambda(|\beta|+)$ denote the limit of $p'_\lambda(x)$ as $x \to |\beta|$ from above, the MM algorithms introduced in the next section are shown to apply to any continuous penalty function $p_\lambda(\beta)$ that is nondecreasing and concave on $(0, \infty)$ such that $p'_\lambda(0+) < \infty$. The previously mentioned penalty functions that satisfy these criteria are hard thresholding, SCAD, LASSO and $L_q$ with $0 < q \leq 1$. Therefore, the methods presented in this article enable a wide range of variable selection algorithms.

Nonetheless, there are some common penalty functions that do not meet our criteria. The entropy penalty is excluded because it is discontinuous, and in fact maximizing the AIC- or BIC-penalized likelihood in cases other than linear regression often requires exhaustive fitting of all possible models. For $q > 1$ the $L_q$ penalty is excluded because it is not concave on $(0, \infty)$; however, the fact that $p_\lambda(|\beta|) = |\beta|^q$ is everywhere differentiable suggests that the penalized likelihood function may be susceptible to gradient-based methods, and hence alternatives such as MM may be of limited value. In particular, the special case of ridge regression ($q = 2$) admits a closed-form maximizer, a fact we allude to following (3.17). But there is a more subtle reason for excluding $L_q$ penalties with $q > 1$. Note that $p'_\lambda(0+) > 0$ for any of our nonexcluded penalty functions. As Fan and Li [7] point out, this fact (which they call *singularity at the origin* because it implies a discontinuous derivative at zero) ensures that the penalized likelihood has the *sparsity* property: the resulting estimator is automatically a thresholding rule that sets small estimated coefficients to zero, thus reducing model complexity. Sparsity is an important property for any penalized likelihood technique that is to be useful in a variable selection setting.

## 3. Maximized penalized likelihood via MM.

It is sometimes a challenging task to find the maximum penalized likelihood estimate. Fan and Li [7] propose a local quadratic approximation for the penalty function: Suppose that we are given an initial value $\boldsymbol{\beta}^{(0)}$. If $\beta_j^{(0)}$ is very close to 0, then set $\hat{\beta}_j = 0$; otherwise, the penalty function is locally approximated by a quadratic function using

$$\frac{\partial}{\partial \beta_j}\{p_\lambda(|\beta_j|)\} = p'_\lambda(|\beta_j|+)\,\mathrm{sgn}(\beta_j) \approx \frac{p'_\lambda(|\beta_j^{(0)}|+)}{|\beta_j^{(0)}|}\beta_j$$

when $\beta_j^{(0)} \neq 0$. In other words,

$$p_\lambda(|\beta_j|) \approx p_\lambda(|\beta_j^{(0)}|) + \frac{\{\beta_j^2 - (\beta_j^{(0)})^2\} p'_\lambda(|\beta_j^{(0)}|+)}{2|\beta_j^{(0)}|} \tag{3.1}$$

for $\beta_j \approx \beta_j^{(0)}$. With the aid of this local quadratic approximation, a Newton–Raphson algorithm (for example) can be used to maximize the penalized likelihood function, where each iteration updates the local quadratic approximation.

In this section, we show that this local quadratic approximation idea is an instance of an MM algorithm. This fact enables us to study the convergence properties of the algorithm using techniques applicable to MM algorithms in general. Throughout this section we refrain from specifying the form of $p_\lambda(\cdot)$, since the derivations apply equally to any one of hard thresholding, LASSO, bridge regression using $L_q$ with $0 < q \leq 1$, SCAD or any other method with penalty function $p_\lambda(\cdot)$ satisfying the conditions of Proposition 3.1.

### 3.1. *Local quadratic approximation as an MM algorithm.*

MM stands for majorize–minimize or minorize–maximize, depending on the context [14]. EM algorithms [6], in which the E-step may be shown to be equivalent to a minorization step, are the most famous examples of MM algorithms, though there are many examples of MM algorithms that involve neither maximum likelihood nor missing data. Heiser [12] and Lange, Hunter and Yang [17] give partial surveys of work in this area. The apparent ambiguity in allowing MM to have two different meanings is harmless, since any maximization problem may be viewed as a minimization problem by changing the sign of the objective function.

Consider the penalty term $-n\sum_j p_\lambda(|\beta_j|)$ of (2.1), ignoring its minus sign for the moment. Mimicking the idea of (3.1), we define the function

$$\Phi_{\theta_0}(\theta) = p_\lambda(|\theta_0|) + \frac{(\theta^2 - \theta_0^2)p'_\lambda(|\theta_0|+)}{2|\theta_0|}. \tag{3.2}$$

We assume that $p_\lambda(\cdot)$ is piecewise differentiable so that $p'_\lambda(|\theta|+)$ exists for all $\theta$. Thus, $\Phi_{\theta_0}(\theta)$ is a well-defined quadratic function of $\theta$ for all real $\theta_0$ except for $\theta_0 = 0$. Section 3.2 remedies the problem that $\Phi_{\theta_0}(\cdot)$ is undefined when $\theta_0 = 0$.

We are interested in penalty functions $p_\lambda(\theta)$ for which

$$\Phi_{\theta_0}(\theta) \geq p_\lambda(|\theta|) \qquad \text{for all } \theta \text{ with equality when } \theta = \theta_0. \tag{3.3}$$

A function $\Phi_{\theta_0}(\theta)$ satisfying condition (3.3) is said to *majorize* $p_\lambda(|\theta|)$ at $\theta_0$. If the direction of the inequality in condition (3.3) were reversed, then $\Phi_{\theta_0}(\theta)$ would be said to *minorize* $p_\lambda(|\theta|)$ at $\theta_0$.

The driving force behind an MM algorithm is the fact that condition (3.3) implies

$$\Phi_{\theta_0}(\theta) - \Phi_{\theta_0}(\theta_0) \geq p_\lambda(|\theta|) - p_\lambda(|\theta_0|),$$

which in turn gives the *descent property*

$$\Phi_{\theta_0}(\theta) < \Phi_{\theta_0}(\theta_0) \quad \text{implies} \quad p_\lambda(|\theta|) < p_\lambda(|\theta_0|). \tag{3.4}$$

In other words, if $\theta_0$ denotes the current iterate, any decrease in the value of $\Phi_{\theta_0}(\theta)$ guarantees a decrease in the value of $p_\lambda(|\theta|)$. If $\theta_k$ denotes the estimate of the parameter at the $k$th iteration, then an iterative minimization algorithm would exploit the descent property by constructing the majorizing function $\Phi_{\theta_k}(\theta)$, then minimizing it to give $\theta_{k+1}$—hence the name "majorize–minimize algorithm." Proposition 3.1 gives sufficient conditions on the penalty function $p_\lambda(\cdot)$ in order that $\Phi_{\theta_0}(\theta)$ majorizes $p_\lambda(|\theta|)$. Several different penalty functions that satisfy these conditions are depicted in Figure 1 along with their majorizing quadratic functions.

PROPOSITION 3.1. *Suppose that on* $(0,\infty)$ $p_\lambda(\cdot)$ *is piecewise differentiable, nondecreasing and concave. Furthermore,* $p_\lambda(\cdot)$ *is continuous at* 0 *and* $p'_\lambda(0+) < \infty$. *Then for all* $\theta_0 \neq 0$, $\Phi_{\theta_0}(\theta)$ *as defined in* (3.2) *majorizes* $p_\lambda(|\theta|)$ *at the points* $\pm|\theta_0|$. *In particular, conditions* (3.3) *and* (3.4) *hold.*

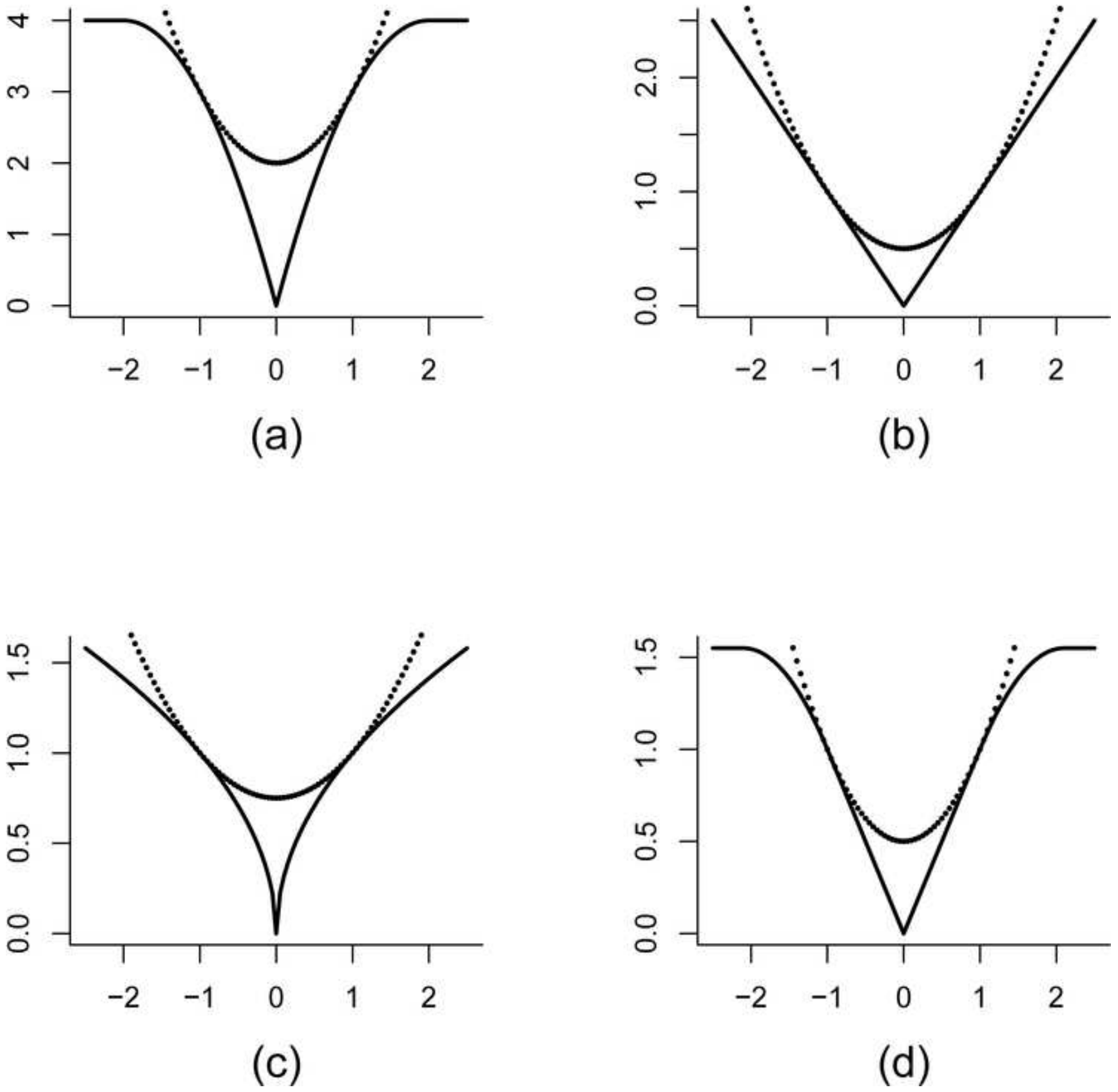

Fig. 1. *Majorizing functions* $\Phi_{\theta_0}(\theta)$ *for various penalty functions are shown as dotted curves; the penalty functions are shown as solid curves. The four penalties are* (a) *hard thresholding with* $\lambda = 2$; (b) $L_1$ *with* $\lambda = 1$; (c) $L_{0.5}$ *with* $\lambda = 1$; *and* (d) *SCAD with* $a = 2.1$ *and* $\lambda = 1$. *In each case* $\theta_0 = 1$.

Next, suppose that we wish to employ the local quadratic approximation idea in an iterative algorithm, where $\boldsymbol{\beta}^{(k)} = (\beta_1^{(k)}, \dots, \beta_d^{(k)})$ denotes the value of $\boldsymbol{\beta}$ at the $k$th iteration. Appending negative signs to $p_\lambda(|\beta_j|)$ and $\Phi_{\beta_j^{(k)}}(\beta_j)$ to turn majorization into minorization, we obtain the following corollary from Proposition 3.1 and (2.1).

Corollary 3.1. *Suppose that* $\beta_j^{(k)} \neq 0$ *for all* $j$ *and that* $p_\lambda(\theta)$ *satisfies the conditions given in Proposition* 3.1. *Then*

$$S_k(\boldsymbol{\beta}) \equiv \ell(\boldsymbol{\beta}) - n \sum_{j=1}^{d} \Phi_{\beta_j^{(k)}}(\beta_j) \tag{3.5}$$

*minorizes* $Q(\boldsymbol{\beta})$ *at* $\boldsymbol{\beta}^{(k)}$.

By the *ascent property*—the analogue for minorizing functions of the descent property (3.4)—Corollary 3.1 suggests that given $\boldsymbol{\beta}^{(k)}$, we should define $\boldsymbol{\beta}^{(k+1)}$ to be the maximizer of $S_k(\boldsymbol{\beta})$, thereby ensuring that $Q(\boldsymbol{\beta}^{(k+1)}) > Q(\boldsymbol{\beta}^{(k)})$. The benefit of replacing one maximization problem by another in this way is that $S_k(\boldsymbol{\beta})$ is susceptible to a gradient-based scheme such as Newton–Raphson, unlike the nondifferentiable function $Q(\boldsymbol{\beta})$. Since the sum in (3.5) is a quadratic function of $\boldsymbol{\beta}$—in fact, the Hessian matrix of this sum is a diagonal matrix—this sum presents no difficulties for maximization. Therefore, the difficulty of maximizing $S_k(\boldsymbol{\beta})$ is determined solely by the form of $\ell(\boldsymbol{\beta})$. For example, in the special case of a linear regression model with normally distributed errors, the log-likelihood function $\ell(\boldsymbol{\beta})$ is itself quadratic, which implies that $S_k(\boldsymbol{\beta})$ may be maximized analytically.

If some of the components of $\boldsymbol{\beta}^{(k)}$ equal zero (or in practice, if some of them are close to zero), the algorithm proceeds by simply setting the final estimates of those components to be zero, deleting them from consideration, then defining the function $S_k(\tilde{\boldsymbol{\beta}})$ as in (3.5), where $\tilde{\boldsymbol{\beta}}$ is the vector composed of the nonzero components of $\boldsymbol{\beta}$. The weakness of this scheme is that once a component is set to zero, it may never reenter the model at a later stage of the algorithm. The modification proposed in Section 3.2 eliminates this weakness.

3.2. *An improved version of local quadratic approximation.* The drawback of $\Phi_{\theta_0}(\theta)$ in (3.2) is that when $\theta_0 = 0$, the denominator $2|\theta_0|$ makes $\Phi_{\theta_0}(\theta)$ undefined. We therefore replace $2|\theta_0|$ by $2(\varepsilon + |\theta_0|)$ for some $\varepsilon > 0$. The resulting perturbed version of $\Phi_{\theta_0}(\theta)$, which is defined for all real $\theta_0$, is no longer a majorizer of $p_\lambda(\theta)$ as required by the MM theory. Nonetheless, we may show that it majorizes a perturbed version of $p_\lambda(\theta)$, which may therefore be used to define a new objective function $Q_\varepsilon(\boldsymbol{\beta})$ that is similar to $Q(\boldsymbol{\beta})$. To this end, we define

$$p_{\lambda,\varepsilon}(|\theta|) = p_\lambda(|\theta|) - \varepsilon \int_0^{|\theta|} \frac{p'_\lambda(t)}{\varepsilon + t}\,dt \tag{3.6}$$

and

$$Q_\varepsilon(\boldsymbol{\beta}) = \ell(\boldsymbol{\beta}) - n \sum_{j=1}^{d} p_{\lambda,\varepsilon}(|\beta_j|). \tag{3.7}$$

The next proposition shows that an MM algorithm may be applied to the maximization of $Q_\varepsilon(\boldsymbol{\beta})$ and suggests that a maximizer of $Q_\varepsilon(\boldsymbol{\beta})$ should be close to a maximizer of $Q(\boldsymbol{\beta})$ as long as $\varepsilon$ is small and $Q(\boldsymbol{\beta})$ is not too flat in the neighborhood of the maximizer.

PROPOSITION 3.2. *Suppose that $p_\lambda(\cdot)$ satisfies the conditions of Proposition 3.1. For $\varepsilon > 0$ define*

$$\Phi_{\theta_0,\varepsilon}(\theta) = p_{\lambda,\varepsilon}(|\theta_0|) + \frac{(\theta^2 - \theta_0^2)p'_\lambda(|\theta_0|+)}{2(\varepsilon + |\theta_0|)}. \tag{3.8}$$

*Then:*

(a) *For any fixed $\varepsilon > 0$,*

$$S_{k,\varepsilon}(\boldsymbol{\beta}) \equiv \ell(\boldsymbol{\beta}) - n \sum_{j=1}^{d} \Phi_{\beta_j^{(k)},\varepsilon}(\beta_j) \tag{3.9}$$

*minorizes $Q_\varepsilon(\boldsymbol{\beta})$ at $\boldsymbol{\beta}^{(k)}$.*

(b) *As $\varepsilon \downarrow 0$, $|Q_\varepsilon(\boldsymbol{\beta}) - Q(\boldsymbol{\beta})| \to 0$ uniformly on compact subsets of the parameter space.*

Note that when the MM algorithm converges and $S_{k,\varepsilon}(\boldsymbol{\beta})$ is maximized by $\boldsymbol{\beta}^{(k)}$, it follows by straightforward differentiation that

$$\frac{\partial S_{k,\varepsilon}(\boldsymbol{\beta}^{(k)})}{\partial \beta_j} = \frac{\partial \ell(\boldsymbol{\beta}^{(k)})}{\partial \beta_j} - np'_\lambda(|\beta_j^{(k)}|)\operatorname{sgn}(\beta_j^{(k)})\frac{|\beta_j^{(k)}|}{\varepsilon + |\beta_j^{(k)}|} = 0.$$

Thus, when $\varepsilon$ is small, the resulting estimator $\hat{\boldsymbol{\beta}}$ approximately satisfies the penalized likelihood equation

$$\frac{\partial Q(\hat{\boldsymbol{\beta}})}{\partial \beta_j} = \frac{\partial \ell(\hat{\boldsymbol{\beta}})}{\partial \beta_j} - n\operatorname{sgn}(\hat{\beta}_j)p'_\lambda(|\hat{\beta}_j|) = 0. \tag{3.10}$$

Suppose that $\hat{\boldsymbol{\beta}}_\varepsilon$ denotes a maximizer of $Q_\varepsilon(\boldsymbol{\beta})$ and $\hat{\boldsymbol{\beta}}_0$ denotes a maximizer of $Q(\boldsymbol{\beta})$. In general, it is impossible to bound $\|\hat{\boldsymbol{\beta}}_\varepsilon - \hat{\boldsymbol{\beta}}_0\|$ as a function of $\varepsilon$ because $Q(\boldsymbol{\beta})$ may be quite flat near its maximum. However, suppose that $Q_\varepsilon(\boldsymbol{\beta})$ is *upper compact*, which means that $\{\boldsymbol{\beta} : Q_\varepsilon(\boldsymbol{\beta}) \geq c\}$ is a compact subset of $R^d$ for any real constant $c$. In this case we may obtain the following corollary of Proposition 3.2(b).

COROLLARY 3.2. *Suppose that $\hat{\boldsymbol{\beta}}_\varepsilon$ denotes a maximizer of $Q_\varepsilon(\boldsymbol{\beta})$. If $Q_\varepsilon(\boldsymbol{\beta})$ is upper compact for all $\varepsilon \geq 0$, then under the conditions of Proposition 3.1 any limit point of the sequence $\{\hat{\boldsymbol{\beta}}_\varepsilon\}_{\varepsilon\downarrow 0}$ is a maximizer of $Q(\boldsymbol{\beta})$.*

Both Proposition 3.2(a) and Corollary 3.2 give results as $\varepsilon \downarrow 0$, which suggests the use of an algorithm in which $\varepsilon$ is allowed to go to zero as the iterations progress. Certainly it would be possible to implement such an algorithm. However, in this article we interpret these results merely as theoretical justification for using the $\varepsilon$ perturbation in the first place, and instead we hold $\varepsilon$ fixed throughout the algorithms we discuss. The choice of this fixed $\varepsilon$ is the subject of Section 3.3.

3.3. *Choice of $\varepsilon$.* Essentially, we want to solve the penalized likelihood equation (3.10) for $\hat{\beta}_j \neq 0$ [recall that $p_\lambda(\beta_j)$ is not differentiable at $\beta_j = 0$]. Suppose, therefore, that convergence is declared in a numerical algorithm whenever $|\partial Q(\hat{\boldsymbol{\beta}})/\partial\beta_j| < \tau$ for a predetermined small tolerance $\tau$. Our algorithm accomplishes this by declaring convergence whenever $|\partial Q_\varepsilon(\hat{\boldsymbol{\beta}})/\partial\beta_j| < \tau/2$, where $\varepsilon$ satisfies

$$\left|\frac{\partial Q_\varepsilon(\hat{\boldsymbol{\beta}})}{\partial\beta_j} - \frac{\partial Q(\hat{\boldsymbol{\beta}})}{\partial\beta_j}\right| < \frac{\tau}{2}. \tag{3.11}$$

Since $p'_\lambda(\theta)$ is nonincreasing on $(0,\infty)$, (3.6) implies

$$\left|\frac{\partial Q_\varepsilon(\hat{\boldsymbol{\beta}})}{\partial\beta_j} - \frac{\partial Q(\hat{\boldsymbol{\beta}})}{\partial\beta_j}\right| = \frac{\varepsilon n p'_\lambda(|\hat{\beta}_j|)}{|\hat{\beta}_j| + \varepsilon} \leq \frac{\varepsilon n p'_\lambda(0+)}{|\hat{\beta}_j|}$$

for $\hat{\boldsymbol{\beta}}_j \neq 0$. Thus, to ensure that (3.11) is satisfied, simply take

$$\varepsilon = \frac{\tau}{2np'_\lambda(0+)} \min\{|\beta_j| : \beta_j \neq 0\}.$$

This may of course lead to a different value of $\varepsilon$ each time $\boldsymbol{\beta}$ changes; therefore, in our implementations we fix

$$\varepsilon = \frac{\tau}{2np'_\lambda(0+)} \min\{|\beta_j^{(0)}| : \beta_j^{(0)} \neq 0\}. \tag{3.12}$$

When the algorithm converges, if $|\partial Q(\hat{\boldsymbol{\beta}})/\partial\beta_j| > \tau$, $\beta_j$ is presumed to be zero. In the numerical examples of Section 4 we take $\tau = 10^{-8}$.

3.4. *The algorithm.* By the ascent property of an MM algorithm, $\boldsymbol{\beta}^{(k+1)}$ should be defined at the $k$th iteration so that

$$S_{k,\varepsilon}(\boldsymbol{\beta}^{(k+1)}) > S_{k,\varepsilon}(\boldsymbol{\beta}^{(k)}). \tag{3.13}$$

Note that if $\boldsymbol{\beta}^{(k+1)}$ satisfies (3.13) without actually maximizing $S_{k,\varepsilon}(\boldsymbol{\beta})$, we still refer to the algorithm as an MM algorithm, even though the second M—for "maximize"—is not quite accurate. Alternatively, we could adopt the convention used for EM algorithms by Dempster, Laird and Rubin [6] and refer to such algorithms as generalized MM, or GMM, algorithms; however, in this article we prefer to require extra duty of the label MM and avoid further crowding the field of acronym-named algorithms.

From (3.9) we see that $S_{k,\varepsilon}(\boldsymbol{\beta})$ consists of two parts, $\ell(\boldsymbol{\beta})$ and the sum of quadratic functions $-\Phi_{\beta_j^{(k)},\varepsilon}(\beta_j)$ of the components of $\boldsymbol{\beta}$. The latter part is easy to maximize directly; thus, the difficulty of maximizing $S_{k,\varepsilon}(\boldsymbol{\beta})$, or at least attaining (3.13), is solely determined by the form of $\ell(\boldsymbol{\beta})$. In general, when $\ell(\boldsymbol{\beta})$ is easy to maximize, then so is $S_{k,\varepsilon}(\boldsymbol{\beta})$, which distinguishes $S_{k,\varepsilon}(\boldsymbol{\beta})$ from the ($\varepsilon$-perturbed) penalized likelihood $Q_\varepsilon(\boldsymbol{\beta})$. Even if $\ell(\boldsymbol{\beta})$ is not easily maximized, at least if it is differentiable, then so is $S_{k,\varepsilon}(\boldsymbol{\beta})$, which means (3.13) may be attained using standard gradient-based maximization methods such as Newton–Raphson. The function $Q_\varepsilon(\boldsymbol{\beta})$, though differentiable, is not easily optimized using gradient-based methods because it is very close to the nondifferentiable function $Q(\boldsymbol{\beta})$.

Although it is impossible to detail all possible forms of likelihood functions $\ell(\boldsymbol{\beta})$ to which these methods apply, we begin with the completely general Newton–Raphson-based algorithm

$$\boldsymbol{\beta}^{(k+1)} = \boldsymbol{\beta}^{(k)} - \alpha_k[\nabla^2 S_{k,\varepsilon}(\boldsymbol{\beta}^{(k)})]^{-1}\nabla S_{k,\varepsilon}(\boldsymbol{\beta}^{(k)}), \tag{3.14}$$

where $\nabla^2 S_{k,\varepsilon}(\cdot)$ and $\nabla S_{k,\varepsilon}(\cdot)$ denote the Hessian matrix and gradient vector, respectively, and $\alpha_k$ is some positive scalar. Using the definition (3.9) of $S_{k,\varepsilon}(\boldsymbol{\beta})$, algorithm (3.14) becomes

$$\boldsymbol{\beta}^{(k+1)} = \boldsymbol{\beta}^{(k)} - \alpha_k[\nabla^2 \ell(\boldsymbol{\beta}^{(k)}) - nE_k]^{-1}[\nabla \ell(\boldsymbol{\beta}^{(k)}) - nE_k\boldsymbol{\beta}^{(k)}], \tag{3.15}$$

where $E_k$ is the diagonal matrix with $(j,j)$th entry $p'_\lambda(|\beta_j^{(k)}|+)/(\varepsilon + |\beta_j^{(k)}|)$. We take the ordinary maximum likelihood estimate to be the initial value $\boldsymbol{\beta}^{(0)}$ in the numerical examples of Section 4.

There are some important special cases. First, we consider perhaps the simplest case but by far the most important case because of its ubiquity—the linear regression model with normal homoscedastic errors, for which

$$\nabla \ell(\boldsymbol{\beta}) = \mathbf{X}^T\mathbf{y} - \mathbf{X}^T\mathbf{X}\boldsymbol{\beta}, \tag{3.16}$$

where $\mathbf{X} = (\mathbf{x}_1, \ldots, \mathbf{x}_n)^T$, the design matrix of the regression model, and $\mathbf{y}$ is the response vector consisting of $y_i$. As pointed out at the beginning of Section 2, we omit mention of the error variance parameter $\sigma^2$ here because this parameter is to be estimated using standard methods once $\boldsymbol{\beta}$ has been estimated. In this case, (3.15) with $\alpha_k = 1$ gives a closed-form maximizer of $S_{k,\varepsilon}(\boldsymbol{\beta})$ because $S_{k,\varepsilon}(\boldsymbol{\beta})$ is exactly a quadratic function. The resulting algorithm

$$\boldsymbol{\beta}^{(k+1)} = \{\mathbf{X}^T\mathbf{X} + nE_k\}^{-1}\mathbf{X}^T\mathbf{y} \tag{3.17}$$

may be viewed as iterative ridge regression. In the case of LASSO, which uses the $L_1$ penalty, this algorithm is guaranteed to converge to the unique maximizer of $Q_\varepsilon(\boldsymbol{\beta})$ (see Corollary 3.3).

The slightly more general case of generalized linear models with canonical link includes common procedures such as logistic regression and Poisson regression. In these cases, the Hessian matrix is

$$\nabla^2 \ell(\boldsymbol{\beta}) = -\mathbf{X}^T V \mathbf{X},$$

where $V$ is a diagonal matrix whose $(i,i)$th entry is given by $v(\mathbf{x}_i^T\boldsymbol{\beta})$ and $v(\cdot)$ is the variance function. Therefore, the Hessian matrix is negative definite provided that $\mathbf{X}$ is of full rank. This means that the vector $-[\nabla^2 S_{k,\varepsilon}(\boldsymbol{\beta}^{(k)})]^{-1}\nabla S_{k,\varepsilon}(\boldsymbol{\beta}^{(k)})$ in (3.14) is a direction of ascent [unless of course $\nabla S_{k,\varepsilon}(\boldsymbol{\beta}^{(k)}) = \mathbf{0}$], which guarantees the existence of a positive $\alpha_k$ such that $\boldsymbol{\beta}^{(k+1)}$ satisfies (3.13). A simple way of determining $\alpha_k$ is the practice of step-halving: try $\alpha_k = 2^{-\nu}$ for $\nu = 0, 1, 2, \ldots$ until the resulting $\boldsymbol{\beta}^{(k+1)}$ satisfies (3.13). For large samples in practice, $\ell(\boldsymbol{\beta})$ tends to be nearly quadratic, particularly in the vicinity of the MLE (which is close to the penalized MLE for large samples), so step-halving does not need to be employed very often. Nonetheless, in our experience it is always wise to check that (3.13) is satisfied; even when the truth of this inequality is guaranteed as in the linear regression model, checking the inequality is often a useful debugging tool. Indeed, whenever possible it is a good idea to check that the objective function itself satisfies $Q_\epsilon(\boldsymbol{\beta}^{(k+1)}) > Q_\epsilon(\boldsymbol{\beta}^{(k)})$; for even though this inequality is guaranteed theoretically by (3.13), in practice many a programming error is caught using this simple check.

In still more general cases, the Hessian matrix $\nabla^2\ell(\boldsymbol{\beta})$ may not be negative definite or it may be difficult to compute. If the Fisher information matrix $I(\boldsymbol{\beta})$ is known, then $-nI(\boldsymbol{\beta})$ may be used in place of $\nabla^2\ell(\boldsymbol{\beta})$. This leads to

$$\boldsymbol{\beta}^{(k+1)} = \boldsymbol{\beta}^{(k)} + \alpha_k[I(\boldsymbol{\beta}^{(k)}) + E_k]^{-1}\left[\frac{1}{n}\nabla\ell(\boldsymbol{\beta}^{(k)}) - E_k\boldsymbol{\beta}^{(k)}\right],$$

and the positive definiteness of the Fisher information will ensure that step-halving will always lead to an increase in $S_{k,\varepsilon}(\boldsymbol{\beta})$.

Finally, we mention the possibility of applying an MM algorithm to a penalized partial likelihood function. Consider the example of the Cox proportional hazards model [4], which is the most popular model in survival data analysis. The variable selection methods of Section 2 are extended to the Cox model by Fan and Li [8]. Let $T_i, C_i$ and $\mathbf{x}_i$ be, respectively, the survival time, the censoring time and the vector of covariates for the $i$th individual. Correspondingly, let $Z_i = \min\{T_i, C_i\}$ be the observed time and let $\delta_i = I(T_i \le C_i)$ be the censoring indicator. It is assumed that $T_i$ and $C_i$ are conditionally independent given $\mathbf{x}_i$ and that the censoring mechanism is noninformative. Under the proportional hazards model, the conditional hazard function $h(t_i|\mathbf{x}_i)$ of $T_i$ given $\mathbf{x}_i$ is given by

$$h(t_i|\mathbf{x}_i) = h_0(t_i)\exp(\mathbf{x}_i^T\boldsymbol{\beta}),$$

where $h_0(t)$ is the baseline hazard function. This is a semiparametric model with parameters $h_0(t)$ and $\boldsymbol{\beta}$. Denote the ordered uncensored failure times by $t_1^0 \le \cdots \le t_N^0$, and let $(j)$ provide the label for the item falling at $t_j^0$, so that the covariates associated with the $N$ failures are $\mathbf{x}_{(1)}, \ldots, \mathbf{x}_{(N)}$. Let $R_j = \{i : Z_i \ge t_j^0\}$ denote the risk set right before the time $t_j^0$. A partial

likelihood is given by

$$\ell_P(\boldsymbol{\beta}) = \sum_{j=1}^{N}\left[\mathbf{x}_{(j)}^T\boldsymbol{\beta} - \log\left\{\sum_{i\in R_j}\exp(\mathbf{x}_i^T\boldsymbol{\beta})\right\}\right].$$

Fan and Li [8] consider variable selection via maximization of the penalized partial likelihood

$$\ell_P(\boldsymbol{\beta}) - n\sum_{j=1}^{d} p_\lambda(|\beta_j|).$$

It can be shown that the Hessian matrix of $\ell_P(\boldsymbol{\beta})$ is negative definite provided that $\mathbf{X}$ is of full rank. Under certain regularity conditions, it can further be shown that in the neighborhood of a maximizer, the partial likelihood is nearly quadratic for large $n$.

3.5. *Convergence.* It is not possible to prove that a generic MM algorithm converges at all, and when an MM algorithm does converge, there is no guarantee that it converges to a global maximum. For example, there are well-known pathological examples in which EM algorithms—or generalized EM algorithms, as discussed following inequality (3.13)—converge to saddle points or fail to converge [18]. Nonetheless, it is often possible to obtain convergence results in specific cases.

We define a stationary point of the function $Q_\varepsilon(\boldsymbol{\beta})$ to be any point $\boldsymbol{\beta}$ at which the gradient vector is zero. Because the differentiable function $S_{k,\varepsilon}(\boldsymbol{\beta})$ is tangent to $Q_\varepsilon(\boldsymbol{\beta})$ at the point $\boldsymbol{\beta}^{(k)}$ by the minorization property, the gradient vectors of $S_{k,\varepsilon}(\boldsymbol{\beta})$ and $Q_\varepsilon(\boldsymbol{\beta})$ are equal when evaluated at $\boldsymbol{\beta}^{(k)}$. Thus, when using the method of Section 3.4 to maximize $S_{k,\varepsilon}(\boldsymbol{\beta})$, we see that fixed points of the algorithm (i.e., points with gradient zero) coincide with stationary points of $Q_\varepsilon(\boldsymbol{\beta})$. Letting $M(\boldsymbol{\beta})$ denote the map implicitly defined by the algorithm that takes $\boldsymbol{\beta}^{(k)}$ to $\boldsymbol{\beta}^{(k+1)}$ for any point $\boldsymbol{\beta}^{(k)}$, inequality (3.13) states that $S_{k,\varepsilon}\{M(\boldsymbol{\beta})\} > S_{k,\varepsilon}(\boldsymbol{\beta})$. The limit points of the set $\{\boldsymbol{\beta}^{(k)} : k = 0, 1, 2, \ldots\}$ are characterized by the following slightly modified version of Lyapunov's theorem [16].

PROPOSITION 3.3. *Given an initial value $\boldsymbol{\beta}^{(0)}$, let $\boldsymbol{\beta}^{(k)} = M^k(\boldsymbol{\beta}^{(0)})$. If $Q_\varepsilon(\boldsymbol{\beta}) = Q_\varepsilon\{M(\boldsymbol{\beta})\}$ only for stationary points $\boldsymbol{\beta}$ of $Q_\varepsilon$ and if $\boldsymbol{\beta}^*$ is a limit point of the sequence $\{\boldsymbol{\beta}^{(k)}\}$ such that $M(\boldsymbol{\beta})$ is continuous at $\boldsymbol{\beta}^*$, then $\boldsymbol{\beta}^*$ is a stationary point of $Q_\varepsilon(\boldsymbol{\beta})$.*

Equation (3.14) with $\alpha_k = 1$ gives

$$M(\boldsymbol{\beta}^{(k)}) = \boldsymbol{\beta}^{(k)} - \{\nabla^2 S_{k,\varepsilon}(\boldsymbol{\beta}^{(k)})\}^{-1}\nabla Q_\varepsilon(\boldsymbol{\beta}^{(k)}), \tag{3.18}$$

where (3.18) uses the fact that $\nabla S_{k,\varepsilon}(\boldsymbol{\beta}^{(k)}) = \nabla Q_\varepsilon(\boldsymbol{\beta}^{(k)})$. As discussed in [16, 17], the derivative of $M(\boldsymbol{\beta})$ gives insight into the local convergence properties of the algorithm. Suppose that $\nabla Q_\varepsilon(\boldsymbol{\beta}^{(k)}) = \mathbf{0}$, so $\boldsymbol{\beta}^{(k)}$ is a stationary point. In this case, differentiating (3.18) gives

$$\nabla M(\boldsymbol{\beta}^{(k)}) = \{\nabla^2 S_{k,\varepsilon}(\boldsymbol{\beta}^{(k)})\}^{-1}\{\nabla^2 S_{k,\varepsilon}(\boldsymbol{\beta}^{(k)}) - \nabla^2 Q_\varepsilon(\boldsymbol{\beta}^{(k)})\}.$$

It is possible to write $\nabla^2 S_{k,\varepsilon}(\boldsymbol{\beta}^{(k)}) - \nabla^2 Q_\varepsilon(\boldsymbol{\beta}^{(k)})$ in closed form as $nA_k$, where $A_k = \operatorname{diag}\{a(\beta_1^{(k)}), \dots, a(\beta_d^{(k)})\}$ and

$$a(t) = \frac{|t|}{\varepsilon + |t|}\left\{p_\lambda''(|t|+) - \frac{p_\lambda(|t|+)}{\varepsilon + |t|}\right\}.$$

Under the conditions of Proposition 3.2, $p_\lambda''(|t|+) \le 0$ and thus $nA_k$ is negative semidefinite, a fact that may also be interpreted as a consequence of the minorization of $Q_\varepsilon(\boldsymbol{\beta})$ by $S_{k,\varepsilon}(\boldsymbol{\beta})$. Furthermore, $\nabla^2 \ell(\boldsymbol{\beta}^{(k)})$ is often negative definite, as pointed out in Section 3.4, which implies that $\nabla^2 S_{k,\varepsilon}(\boldsymbol{\beta}^{(k)})$ is negative definite. This fact, together with the fact that $\nabla^2 S_{k,\varepsilon}(\boldsymbol{\beta}^{(k)}) - \nabla^2 Q_\varepsilon(\boldsymbol{\beta}^{(k)})$ is negative semidefinite, implies that the eigenvalues of $\nabla M(\boldsymbol{\beta}^{(k)})$ are all contained in the interval $[0, 1)$ [13]. Ostrowski's theorem [22] thus implies that the MM algorithm defined by (3.18) is locally attracted to $\boldsymbol{\beta}^{(k)}$ and that the rate of convergence to $\boldsymbol{\beta}^{(k)}$ in a neighborhood of $\boldsymbol{\beta}^{(k)}$ is linear with rate equal to the largest eigenvalue of $\dot{M}(\boldsymbol{\beta}^{(k)})$. In other words, if $\boldsymbol{\beta}^*$ is a stationary point and $\rho < 1$ is the largest eigenvalue of $\nabla M(\boldsymbol{\beta}^*)$, then for any $\delta > 0$ such that $0 < \rho - \delta < \rho + \delta < 1$, there exists a neighborhood $N_\delta$ containing $\boldsymbol{\beta}^*$ such that for all $\boldsymbol{\beta} \in N_\delta$,

$$(\rho - \delta)\|\beta - \beta^*\| \le \|M(\boldsymbol{\beta}) - \beta^*\| \le (\rho + \delta)\|\beta - \beta^*\|.$$

Further details about the rate of convergence for similar algorithms may be found in [16, 17].

Lyapunov's theorem (Proposition 3.3) gives a necessary condition for a point to be a limit point of a particular MM algorithm. To conclude this section, we consider a sufficient condition for the existence of a limit point. Suppose that the function $Q_\varepsilon(\boldsymbol{\beta})$ is upper compact, as defined in Section 3.2. Then given the initial parameter vector $\boldsymbol{\beta}^{(0)}$, the set $B = \{\boldsymbol{\beta} \in R^d : Q_\varepsilon(\boldsymbol{\beta}) \ge Q(\boldsymbol{\beta}^{(0)})\}$ is compact; furthermore, by (3.6) and (3.7), $Q_\varepsilon(\boldsymbol{\beta}) \ge Q(\boldsymbol{\beta})$ so that $B$ contains the entire sequence $\{\boldsymbol{\beta}^{(k)}\}_{k=0}^\infty$. This guarantees that the sequence has at least one limit point, which must therefore be a stationary point of $Q_\varepsilon(\boldsymbol{\beta})$ by Proposition 3.3. If in addition there is no more than one stationary point—for example, if $Q_\varepsilon(\boldsymbol{\beta})$ is strictly concave—then we may conclude that the algorithm must converge to the unique stationary point.

Upper compactness of $Q_\varepsilon(\boldsymbol{\beta})$ follows as long as $Q_\varepsilon(\boldsymbol{\beta}) \to -\infty$ whenever $\|\boldsymbol{\beta}\| \to \infty$; this is often not difficult to prove for specific examples. In the particular case of the $L_1$ penalty (LASSO), strict concavity also holds as long as $\ell(\boldsymbol{\beta})$ is strictly concave, which implies the following corollary.

COROLLARY 3.3. *If $p_\lambda(|\theta|) = \lambda|\theta|$ and $\ell(\boldsymbol{\beta})$ is strictly concave and upper compact, then the MM algorithm of* (3.15) *gives a sequence $\{\boldsymbol{\beta}^{(k)}\}$ converging to the unique maximizer of $Q_\varepsilon(\boldsymbol{\beta})$ for any $\varepsilon > 0$.*

In particular, Corollary 3.3 implies that using our algorithm with the $\varepsilon$-perturbed LASSO penalty guarantees convergence to the maximum penalized likelihood estimator for any full-rank generalized linear model or (say) Cox proportional hazards model. However, strict concavity of $Q_\varepsilon(\boldsymbol{\beta})$ is not typical for other penalty functions presented in this article in light of the requirement in Proposition 3.2 that $p_\lambda(\cdot)$ be concave—and hence that $-p_\lambda(\cdot)$ be convex—on $(0, \infty)$. This fact means that when an MM algorithm using some penalty function other than $L_1$ converges, then it may converge to a local, rather than a global, maximizer of $Q_\varepsilon(\boldsymbol{\beta})$. This can actually be an advantage, since one might like to know if the penalized likelihood has multiple local maxima.

## 4. Numerical examples.

Since Fan and Li [7] have already compared the performance of LASSO with SCAD using a local quadratic approximation and other existing methods, in the following four numerical examples we focus on assessing the performance of the proposed algorithms using the SCAD penalty (2.2). Namely, we compare the unmodified LQA to our modified version (both using SCAD), where $\varepsilon$ is chosen according to (3.12) with $\tau = 10^{-8}$. For SCAD, $p'_\lambda(0+) = \lambda$, and this tuning parameter $\lambda$ is chosen using generalized cross-validation, or GCV [5]. As suggested by Fan and Li [7], we take $a = 3.7$ in the definition of SCAD.

The Newton–Raphson algorithm (3.15) enables a standard error estimate via a sandwich formula:

$$\widehat{\mathrm{cov}}(\hat{\boldsymbol{\beta}}) = \{\nabla^2\ell(\hat{\boldsymbol{\beta}}) - nE_k\}^{-1}\widehat{\mathrm{cov}}\{\nabla S_{k,\varepsilon}(\hat{\boldsymbol{\beta}})\}\{\nabla^2\ell(\hat{\boldsymbol{\beta}}) - nE_k\}^{-1}, \tag{4.1}$$

where

$$\begin{aligned}\widehat{\mathrm{cov}}\{\nabla S_{k,\varepsilon}(\boldsymbol{\beta})\} &= \frac{1}{n}\sum_{i=1}^n[\nabla\ell_i(\boldsymbol{\beta}) - nE_k\boldsymbol{\beta}][\nabla\ell_i(\boldsymbol{\beta}) - nE_k\boldsymbol{\beta}]^T \\ &\qquad - \left[\frac{1}{n}\sum_{i=1}^n\nabla\ell_i(\boldsymbol{\beta}) - nE_k\boldsymbol{\beta}\right]\left[\frac{1}{n}\sum_{i=1}^n\nabla\ell_i(\boldsymbol{\beta}) - nE_k\boldsymbol{\beta}\right]^T \\ &= \frac{1}{n}\sum_{i=1}^n[\nabla\ell_i(\boldsymbol{\beta})][\nabla\ell_i(\boldsymbol{\beta})]^T - \left[\frac{1}{n}\sum_{i=1}^n\nabla\ell_i(\boldsymbol{\beta})\right]\left[\frac{1}{n}\sum_{i=1}^n\nabla\ell_i(\boldsymbol{\beta})\right]^T.\end{aligned}$$

Naturally, another estimate may be formed if $-nI(\hat{\boldsymbol{\beta}})$ is substituted for $\nabla\ell(\hat{\boldsymbol{\beta}})$ in (4.1). Fan and Peng [10] establish the consistency of this sandwich formula for related problems, and their method of proof may be adapted to this situation, though we do not do so in this article.

For the simulated examples, Examples 1–3, we compare the performance of the proposed procedures along with AIC and BIC in terms of model error and model complexity. With $\mu(\mathbf{x}) = E(Y|\mathbf{x})$, model error (ME) is defined as $E\{\hat{\mu}(\mathbf{x}) - \mu(\mathbf{x})\}^2$, where the expectation is taken with respect to a new observation $\mathbf{x}$. The MEs of the underlying procedures are divided by that of the ordinary maximum likelihood estimate, so we report relative model error (RME).

EXAMPLE 1 (*Linear regression*). In this example we generated 500 data sets, each of which consists of 100 observations from the model

$$y = \mathbf{x}^T\boldsymbol{\beta} + \varepsilon,$$

where $\boldsymbol{\beta}$ is a twelve-dimensional vector whose first, fifth and ninth components are 3, 1.5 and 2, respectively, and whose other components equal 0. The components of $\mathbf{x}$ and $\varepsilon$ are standard normal and the correlation between $x_i$ and $x_j$ is taken to be $\rho$. In our simulation we consider three cases: $\rho = 0.1$, $\rho = 0.5$ and $\rho = 0.9$. In this case there is a closed form for the model error, namely $\mathrm{ME}(\hat{\boldsymbol{\beta}}) = (\hat{\boldsymbol{\beta}} - \boldsymbol{\beta})^T \mathrm{cov}(\mathbf{x})(\hat{\boldsymbol{\beta}} - \boldsymbol{\beta})$. The median of the relative model error (RME) over 500 simulated data sets is summarized in Table 1. The average number of 0 coefficients is also reported in Table 1, in which the column labeled "C" gives the average number of coefficients, of the nine true zeros, correctly set to zero and the column labeled "I" gives the average number of the three true nonzeros incorrectly set to zero.

In Table 1 New and LQA refer to the newly proposed algorithm and the local quadratic approximation algorithm of Fan and Li [7]. AIC and BIC stand for the best subset variable selection procedures that minimize AIC scores and BIC scores. Finally, "Oracle" stands for the oracle estimate computed by using the true model $y = \beta_1 x_1 + \beta_5 x_5 + \beta_9 x_9 + \varepsilon$. When the correlation among the covariates is small or moderate, we see that the new algorithm performs the best in terms of model error and LQA also performs very well; their RMEs are both very close to those of the oracle estimator. When the covariates are highly correlated, the new algorithm outperforms LQA in terms of both model error and model complexity. The performance of BIC and AIC remains almost the same for the three cases in this example; Table 1 indicates that BIC performs better than AIC.

TABLE 1
*Relative model errors for linear regression*

| | RME | Zeros | | RME | Zeros | | RME | Zeros | |
|---|---|---|---|---|---|---|---|---|---|
| | **Median** | **C** | **I** | **Median** | **C** | **I** | **Median** | **C** | **I** |
| **Method** | $\boldsymbol{\rho = 0.9}$ | | | $\boldsymbol{\rho = 0.5}$ | | | $\boldsymbol{\rho = 0.1}$ | | |
| New | 0.437 | 8.346 | 0 | 0.215 | 8.708 | 0 | 0.238 | 8.292 | 0 |
| LQA | 0.590 | 7.772 | 0 | 0.237 | 8.680 | 0 | 0.269 | 8.272 | 0 |
| BIC | 0.337 | 8.644 | 0 | 0.335 | 8.652 | 0 | 0.328 | 8.656 | 0 |
| AIC | 0.672 | 7.358 | 0 | 0.673 | 7.324 | 0 | 0.668 | 7.374 | 0 |
| Oracle | 0.201 | 9.000 | 0 | 0.202 | 9 | 0 | 0.211 | 9 | 0 |

TABLE 2
*Standard deviations and standard errors of $\hat{\beta}_1$ in the linear regression model*

| | **SD** | **SE (std(SE))** | **SD** | **SE (std(SE))** | **SD** | **SE (std(SE))** |
|---|---|---|---|---|---|---|
| **Method** | $\rho = 0.9$ | | $\rho = 0.5$ | | $\rho = 0.1$ | |
| LSE | 0.339 | 0.322 (0.035) | 0.152 | 0.144 (0.016) | 0.114 | 0.109 (0.012) |
| New | 0.303 | 0.260 (0.036) | 0.129 | 0.120 (0.017) | 0.104 | 0.098 (0.014) |
| LQA | 0.315 | 0.265 (0.037) | 0.128 | 0.120 (0.017) | 0.105 | 0.098 (0.014) |
| BIC | 0.295 | 0.269 (0.028) | 0.133 | 0.124 (0.013) | 0.105 | 0.101 (0.010) |
| AIC | 0.322 | 0.278 (0.029) | 0.145 | 0.128 (0.013) | 0.109 | 0.101 (0.010) |
| Oracle | 0.270 | 0.264 (0.027) | 0.126 | 0.124 (0.013) | 0.103 | 0.102 (0.010) |

We now test the accuracy of the proposed standard error formula. The standard deviation of the estimated coefficients for the 500 simulated data sets, denoted by SD, can be regarded as the true standard deviation except for Monte Carlo error. The average of the estimated standard errors for the 500 simulated data sets, denoted by SE, and their standard deviation, denoted by std(SE), gauge the overall performance of the standard error formula. Table 2 only presents the SD, SE and std(SE) of $\beta_1$. The results for other coefficients are similar. In Table 2, LSE stands for the ordinary least squares estimate; other notation is the same as that in Table 1. The differences between SD and SE are less than twice std(SE), which suggests that the proposed standard error formula works fairly well. However, the SE appears to consistently underestimate the SD, a common phenomenon (see [15]), so it may benefit from some slight modification.

EXAMPLE 2 (*Logistic regression*). In this example we assess the performance of the proposed algorithm for a logistic regression model. We generated 500 data sets, each of which consists of 200 observations, from the logistic regression model

$$\mu(\mathbf{x}) \equiv \{P(Y = 1|\mathbf{x})\} = \frac{\exp(\mathbf{x}^T \boldsymbol{\beta})}{1 + \exp(\mathbf{x}^T \boldsymbol{\beta})}, \tag{4.2}$$

where $\boldsymbol{\beta}$ is a nine-dimensional vector whose first, fourth and seventh components are 3, 1.5 and 2 respectively, and whose other components equal 0. The components of $\mathbf{x}$ are standard normal, where the correlation between $x_i$ and $x_j$ is $\rho$. In our simulation we consider two cases, $\rho = 0.25$ and $\rho = 0.75$. Unlike the model error for linear regression models, there is no closed form of ME for the logistic regression model in this example. The ME, summarized in Table 3, is estimated by 50,000 Monte Carlo simulations. Notation in Table 3 is the same as that in Table 1. It can be seen from Table 3 that the newly proposed algorithm performs better than LQA in terms of model error and model complexity. We further test the accuracy of the standard error formula derived by using the sandwich formula (4.1). The results are similar to those in Table 2—the proposed standard error formula works fairly well—so they are omitted here.

Best variable subset selection using the BIC criterion is seen in Table 3 to perform quite well relative to other methods. However, best subset selection in this example requires an exhaustive search over all possible subsets, and therefore it is computationally expensive. The methods we propose can dramatically reduce computational cost. To demonstrate this point, random samples of size 200 were generated from model (4.2) with $\boldsymbol{\beta}$ being a $d$-dimensional vector whose first, fourth and seventh components are 3, 1.5 and 2, respectively, and whose other components equal 0. Table 4 depicts the average computing time for each simulation with $d = 8, \ldots, 11$, and indicates that computing times for the BIC and AIC criteria increase exponentially with the dimension $d$, making these methods impractical for parameter sets much larger than those tested here. Given the increasing importance of variable selection problems in fields like genetics and data mining, where the number of variables is measured in the hundreds or even thousands, efficiency of algorithms is an important consideration.

TABLE 3
*Relative model errors for logistic regression*

| | **RME** | **Zeros** | | **RME** | **Zeros** | |
|---|---|---|---|---|---|---|
| | **Median** | **C** | **I** | **Median** | **C** | **I** |
| **Method** | $\boldsymbol{\rho = 0.25}$ | | | $\boldsymbol{\rho = 0.75}$ | | |
| New | 0.277 | 5.922 | 0 | 0.528 | 5.534 | 0.222 |
| LQA | 0.368 | 5.728 | 0 | 0.644 | 4.970 | 0.090 |
| BIC | 0.304 | 5.860 | 0 | 0.399 | 5.796 | 0.304 |
| AIC | 0.673 | 4.930 | 0 | 0.683 | 4.860 | 0.092 |
| Oracle | 0.241 | 6 | 0 | 0.216 | 6 | 0 |

TABLE 4
*Computing time for the logistic model (seconds per simulation)*

| $\rho$ | $d$ | **New** | **LQA** | **BIC** | **AIC** |
|---|---|---|---|---|---|
| 0.25 | 8 | 0.287 | 0.142 | 2.701 | 2.699 |
| | 9 | 0.316 | 0.151 | 5.702 | 5.694 |
| | 10 | 0.348 | 0.180 | 11.761 | 11.754 |
| | 11 | 0.424 | 0.199 | 26.702 | 26.576 |
| 0.75 | 8 | 0.395 | 0.199 | 2.171 | 2.166 |
| | 9 | 0.438 | 0.205 | 4.554 | 4.553 |
| | 10 | 0.452 | 0.225 | 9.499 | 9.518 |
| | 11 | 0.532 | 0.244 | 19.915 | 19.959 |

EXAMPLE 3 (*Cox model*). We investigate the performance of the proposed algorithm for the Cox proportional hazard model in this example. We simulated 500 data sets each for sample sizes $n = 40, 50$ and 60 from the exponential hazard model

$$h(t|\mathbf{x}) = \exp(\mathbf{x}^T\boldsymbol{\beta}), \tag{4.3}$$

where $\boldsymbol{\beta} = (0.8, 0, 0, 1, 0, 0, 0.6, 0)^T$. This model is used in the simulation study of Fan and Li [8]. The $x_u$'s were marginally standard normal and the correlation between $x_u$ and $x_v$ was $\rho^{|u-v|}$ with $\rho = 0.5$. The distribution of the censoring time is an exponential distribution with mean $U\exp(\mathbf{x}^T\boldsymbol{\beta}_0)$, where $U$ is randomly generated from the uniform distribution over $[1, 3]$ for each simulated data set, so that 30% of the data are censored. Here $\boldsymbol{\beta}_0 = \boldsymbol{\beta}$ is regarded as a known constant so that the censoring scheme is noninformative. The model error $E\{\hat{\mu}(\mathbf{x}) - \mu(\mathbf{x})\}^2$ is estimated by 50,000 Monte Carlo simulations and is summarized in Table 5. The performance of the newly proposed algorithm is similar to that of LQA. Both the new algorithm and LQA perform better than best subset variable selection with the AIC or BIC criteria. Note that the BIC criterion is a consistent variable selection criterion. Therefore, as the sample size increases, its performance becomes closer to that of the nonconcave penalized partial likelihood procedures. We also test the accuracy of the standard error formula derived by using the sandwich formula (4.1). The results are similar to those in Table 2; the proposed standard error formula works fairly well.

As in Example 2, we take $\boldsymbol{\beta}$ to be a $d$-dimensional vector whose first, fourth and seventh components are nonzero (0.8, 1.0, and 0.6, resp.) and whose other components equal 0. Table 6 shows that the proposed algorithm and LQA can dramatically save computing time compared with AIC and BIC.

As a referee pointed out, it is of interest to investigate the performance of variable selection algorithms when the "full" model is misspecified. Model misspecification is a concern for all variable selection procedures, not merely those discussed in this article. To address this issue, we generated a random sample from model (4.3) with $\boldsymbol{\beta} = (0.8, 0, 0, 1, 0, 0, 0.6, 0, \beta_9, \beta_{10})^T$, where $\beta_9 = \beta_{10} = 0.2$ or 0.4. The first eight components of $\mathbf{x}$ are the same as those in the above simulation. We take $x_9 = (x_1^2 - 1)/\sqrt{2}$ and $x_{10} = (x_2^2 - 1)/\sqrt{2}$. In our model fitting, our "full" model only uses the first eight components of $\mathbf{x}$. Thus, we misspecified the full model by ignoring the last two components. Based on the "full" model, variable selection procedures are carried out. The oracle procedure uses $(x_1, x_4, x_7, x_9, x_{10})^T$ to fit the data. The model error $E\{\hat{\mu}(\mathbf{x}) - \mu(\mathbf{x})\}^2$, where $\mu(\mathbf{x})$ is the mean function of the true model including all ten components of $\mathbf{x}$, is estimated by 50,000 Monte Carlo simulations and is summarized in Table 7, from which we can see that all variable selection procedures outperform the full model. This implies that selecting significant variables can dramatically reduce both model error and model complexity. From Table 7 we can see that both the newly proposed algorithm and LQA significantly reduce the model error of best subset variable selection using AIC or BIC.

Table 5
*Relative model errors for the Cox model*

| | RME | Zeros | | RME | Zeros | | RME | Zeros | |
|---|---|---|---|---|---|---|---|---|---|
| | Median | C | I | Median | C | I | Median | C | I |
| **Method** | $n = 40$ | | | $n = 50$ | | | $n = 60$ | | |
| New | 0.173 | 4.790 | 1.396 | 0.288 | 4.818 | 1.144 | 0.324 | 4.882 | 0.904 |
| LQA | 0.174 | 4.260 | 0.626 | 0.296 | 4.288 | 0.440 | 0.303 | 4.332 | 0.260 |
| BIC | 0.247 | 4.492 | 0.606 | 0.337 | 4.564 | 0.442 | 0.344 | 4.624 | 0.272 |
| AIC | 0.470 | 3.880 | 0.358 | 0.551 | 3.948 | 0.240 | 0.577 | 3.986 | 0.160 |
| Oracle | 0.103 | 5 | 0 | 0.152 | 5 | 0 | 0.215 | 5 | 0 |

Table 6
*Computing time for the Cox model (seconds per simulation)*

| $n$ | $d$ | New | LQA | BIC | AIC |
|---|---|---|---|---|---|
| 40 | 8 | 0.248 | 0.147 | 0.415 | 0.418 |
| | 9 | 0.258 | 0.140 | 0.843 | 0.843 |
| | 10 | 0.299 | 0.149 | 1.711 | 1.712 |
| | 11 | 0.327 | 0.162 | 3.680 | 3.675 |
| 50 | 8 | 0.320 | 0.197 | 0.588 | 0.591 |
| | 9 | 0.364 | 0.200 | 1.218 | 1.225 |
| | 10 | 0.406 | 0.217 | 2.532 | 2.531 |
| | 11 | 0.466 | 0.211 | 5.189 | 5.171 |
| 60 | 8 | 0.417 | 0.263 | 0.820 | 0.827 |
| | 9 | 0.474 | 0.268 | 1.795 | 1.768 |
| | 10 | 0.513 | 0.279 | 3.454 | 3.454 |
| | 11 | 0.574 | 0.288 | 7.219 | 7.162 |

EXAMPLE 4 (*Environmental data*). In this example, we illustrate the proposed algorithm in the context of analysis of an environmental data set. This data set consists of the number of daily hospital admissions for circulation and respiration problems and daily measurements of air pollutants. It was collected in Hong Kong from January 1, 1994 to December 31, 1995 (courtesy of T. S. Lau). Of interest is the association between levels of pollutants and the total number of daily hospital admissions for circulatory and respiratory problems. The response is the number of admissions and the covariates $X_1$ to $X_3$ are the levels (in $\mu$g/m$^3$) of the pollutants sulfur dioxide, nitrogen dioxide and dust. Because the response is count data, it is reasonable to use a Poisson regression model with mean $\mu(\mathbf{x})$ to analyze this data set. To reduce modeling bias, we include all linear, quadratic and interaction terms among the three air pollutants in our full model. Since empirical studies show that there may be a trend over time, we allow an intercept depending on time, the date on which observations were collected. In other words, we consider the model

$$\begin{aligned}\log\{\mu(\mathbf{x})\} = \beta_0(t) + \beta_1 X_1 + \beta_2 X_2 + \beta_3 X_3 + \beta_4 X_1^2 + \beta_5 X_2^2 \\ + \beta_6 X_3^2 + \beta_7 X_1 X_2 + \beta_8 X_1 X_3 + \beta_9 X_2 X_3.\end{aligned}$$

TABLE 7

*Relative model errors for misspecified Cox models*

| | | **RME** | **# of zeros** | **RME** | **# of zeros** | **RME** | **# of zeros** |
|---|---|---|---|---|---|---|---|
| $(\boldsymbol{\beta_9}, \boldsymbol{\beta_{10}})$ | **Method** | $\boldsymbol{n = 40}$ | | $\boldsymbol{n = 50}$ | | $\boldsymbol{n = 60}$ | |
| (0.2, 0.2) | New | 0.155 | 6.276 | 0.259 | 6.142 | 0.298 | 5.906 |
| | LQA | 0.146 | 4.992 | 0.251 | 4.822 | 0.294 | 4.662 |
| | BIC | 0.244 | 5.186 | 0.356 | 5.020 | 0.487 | 4.932 |
| | AIC | 0.441 | 4.336 | 0.564 | 4.162 | 0.657 | 4.104 |
| | Oracle | 0.254 | 5 | 0.194 | 5 | 0.239 | 5 |
| (0.4, 0.4) | New | 0.205 | 6.544 | 0.323 | 6.384 | 0.479 | 6.132 |
| | LQA | 0.197 | 5.166 | 0.324 | 4.924 | 0.469 | 4.806 |
| | BIC | 0.345 | 5.352 | 0.477 | 5.156 | 0.593 | 5.110 |
| | AIC | 0.560 | 4.444 | 0.670 | 4.236 | 0.755 | 4.240 |
| | Oracle | 0.268 | 5 | 0.228 | 5 | 0.273 | 5 |

We further parameterize the intercept function by a cubic spline

$$\beta_0(t) = \beta_{00} + \beta_{01}t + \beta_{02}t^2 + \beta_{03}t^3 + \sum_{j=1}^{5} \beta_{0(j+4)}(t - k_j)_+^3,$$

where the knots $k_j$ are chosen to be the 10th, 25th, 50th, 75th and 90th percentiles of $t$. Thus, we are dealing with a Poisson regression model with 18 variables. To avoid numerical instability, the time variable and the air pollutant variables are standardized.

Generalized cross-validation is used to select the tuning parameter $\lambda$ for the new algorithm using SCAD. The plot of the GCV scores against $\lambda$ is depicted in Figure 2(a), and the selected $\lambda$ equals 0.1933. In Table 8, we see that all linear terms are very significant, whereas the quadratic terms of $SO_2$ and dust and the $SO_2 \times$ dust interaction are not significant. The plot of the estimated intercept function $\beta_0(t)$ is depicted in Figure 2(b) along with the estimated intercept function under the full model. The two estimated intercept functions are almost identical and capture the time trend very well, but the new algorithm saves two degrees of freedom by deleting the $t^3$ and $(t - k_3)_+^3$ terms from the intercept function.

For the purpose of comparison, SCAD using the LQA is also applied to this data set. The plot of the GCV scores against $\lambda$ is also depicted in Figure 2(a), and the selected $\lambda$ again equals 0.1933. In this case, not only the $t^3$ and

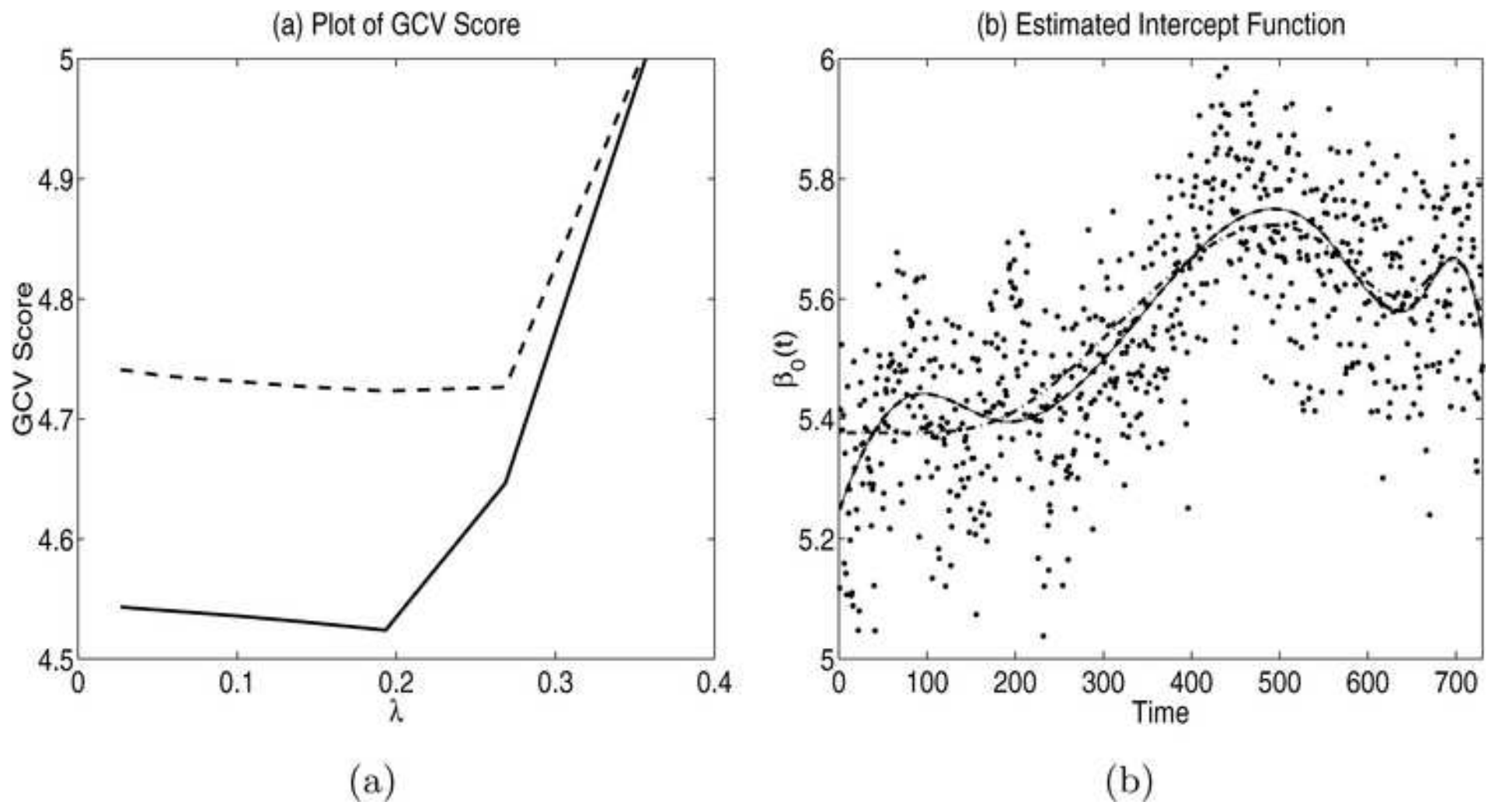

Fig. 2. *In* (a) *the solid line indicates the GCV scores for SCAD using the new algorithm, and the dashed line indicates the same thing for the LQA algorithm. In* (b) *the solid and thick dashed lines that nearly coincide indicate the estimated intercept functions for the new algorithm and the full model, respectively; the dash-dotted line is for LQA. The dots are the full-model residuals* $\log(y) - \mathbf{x}^T\hat{\boldsymbol{\beta}}_{\text{MLE}}$.

TABLE 8
*Estimated coefficients and their standard errors*

| **Covariate** | **MLE** | **New** | **LQA** |
|---|---|---|---|
| $SO_2$ | 0.0082 (0.0041) | 0.0043 (0.0024) | 0.0090 (0.0029) |
| $NO_2$ | 0.0238 (0.0051) | 0.0260 (0.0037) | 0.0311 (0.0033) |
| Dust | 0.0195 (0.0054) | 0.0173 (0.0037) | 0.0043 (0.0026) |
| $SO_2^2$ | −0.0029 (0.0013) | 0 (0.0009) | −0.0025 (0.0010) |
| $NO_2^2$ | 0.0204 (0.0043) | 0.0118 (0.0029) | 0.0157 (0.0031) |
| Dust$^2$ | 0.0042 (0.0025) | 0 (0.0015) | 0.0060 (0.0018) |
| $SO_2 \times NO_2$ | −0.0120 (0.0039) | −0.0050 (0.0021) | −0.0074 (0.0022) |
| $SO_2 \times$ dust | 0.0086 (0.0047) | 0 (<0.00005) | 0 (NA) |
| $NO_2 \times$ dust | −0.0305 (0.0061) | −0.0176 (0.0037) | −0.0262 (0.0041) |

NOTE: NA stands for "Not available."

$(t-k_3)_+^3$ terms but also the $t$ and $t^2$ terms are deleted from the intercept function. The estimated intercept function is the dash-dotted curve in Figure 2(b), and now the resulting estimated intercept function looks dramatically different from the one estimated under the full model, and furthermore it appears to do a poor job of capturing the overall trend. The LQA estimates shown in Table 8 are quite different from those obtained using the new algorithm, even though they both use the same tuning parameter. Recall that LQA suffers the drawback that once a parameter is deleted, it cannot reenter the model, which appears to have had a major impact on the LQA model estimates in this case. Standard errors in Table 8 are available for the deleted coefficients in the new model, but not LQA because the two algorithms use different deletion criteria.

**5. Discussion.** We have shown how a particular class of MM algorithms may be applied to variable selection. In modifying previous work on variable selection using penalized least squares and penalized likelihood by Fan and Li [7, 8, 9] and Cai, Fan, Li and Zhou [3], we have shown how a slight perturbation of the penalty function can eliminate the possibility of mistakenly excluding variables too soon while simultaneously enabling certain convergence results. While the numerical tests given here deal with four very diverse models, the range of possible applications of this method is even broader. Generally speaking, the MM algorithms of this article may be applied to any situation where an objective function—whether a likelihood or not—is penalized using a penalty function such as the one in (2.1). If the goal is to maximize the penalized objective function and $p_\lambda(\cdot)$ satisfies the conditions of Proposition 3.1, then an MM algorithm may be applicable. If the original (unpenalized) objective function is concave, then the modified Newton–Raphson approach of Section 3.4 holds promise. In Section 2 we list several distinct classes of penalty functions in the literature that satisfy the conditions of Proposition 3.1, but there may also be other useful penalties to which our method applies.

The numerical tests of Section 4 indicate that the modified SCAD penalty we propose performs well on simulated data sets. This algorithm has comparable relative model error (sometimes quite a bit better, sometimes slightly worse) to BIC and the unmodified SCAD penalty implemented using LQA, and it consistently outperforms AIC. Our proposed algorithm tends to result in more parsimonious models than the LQA algorithm, typically identifying more actual zeros correctly but also eliminating too many nonzero coefficients. This fact is surprising in light of the drawback that LQA can exclude variables too soon during the iterative process, a drawback that our algorithm corrects. The particular choice of $\varepsilon$, addressed in Section 3.3, may warrant further study because of its influence on the complexity of the final model chosen.

An important difference between our algorithm and both AIC and BIC is the fact that the latter two methods are often not computationally efficient, with computing time scaling exponentially in the number of candidate variables whenever it becomes necessary to search exhaustively over the whole model space. This means that in problems with hundreds or thousands of candidate variables, AIC and BIC can be difficult if not impossible to implement. Such problems are becoming more and more prevalent in the statistical literature as topics such as microarray data and data mining increase in popularity.

Finally, we have seen in one example involving the Cox proportional hazards model that both our method and LQA perform well when the model is misspecified, even outperforming the oracle method for samples of size 40. Although questions of model misspecification are largely outside the scope of this article, it is useful to remember that although simulation studies can aid our understanding, the model assumed for any real data set is only an approximation of reality.

## APPENDIX

PROOF OF PROPOSITION 3.1. The proof uses the following lemma.

LEMMA A.1. *Under the assumptions of Proposition 3.1, $p'_\lambda(\theta+)/(\varepsilon+\theta)$ is a nonincreasing function of $\theta$ for any nonnegative $\varepsilon$.*

The proof of the lemma is immediate: Both $p'_\lambda(\theta)$ and $(\varepsilon+\theta)^{-1}$ are positive and nonincreasing on $(0,\infty)$, so their product is nonincreasing. For any $\theta > 0$, we see that

$$\lim_{x\to\theta+} \frac{d}{dx}[\Phi_{\theta_0}(x) - p_\lambda(|x|)] = \theta\left[\frac{p'_\lambda(|\theta_0|+)}{|\theta_0|} - \frac{p'_\lambda(\theta+)}{\theta}\right]. \tag{A.1}$$

Furthermore, since $p_\lambda(\cdot)$ is nondecreasing and concave on $(0,\infty)$ (and continuous at 0), it is also continuous on $[0,\infty)$. Thus, $\Phi_{\theta_0}(\theta) - p_\lambda(|\theta|)$ is an even function, piecewise differentiable and continuous everywhere. Taking $\varepsilon = 0$ in Lemma A.1, (A.1) implies that $\Phi_{\theta_0}(\theta) - p_\lambda(|\theta|)$ is nonincreasing for $\theta \in (0, |\theta_0|)$ and nondecreasing for $\theta \in (|\theta_0|, \infty)$; this function is therefore minimized on $(0,\infty)$ at $|\theta_0|$. Since it is clear that $\Phi_{\theta_0}(|\theta_0|) = p_\lambda(|\theta_0|)$ and $\Phi_{\theta_0}(-|\theta_0|) = p_\lambda(-|\theta_0|)$, condition (3.3) is satisfied for $\theta_0 = \pm|\theta_0|$. □

Proof of Proposition 3.2. For part (a), it suffices to show that $\Phi_{\theta_0,\varepsilon}(\theta)$ majorizes $p_{\lambda,\varepsilon}(|\theta|)$ at $\theta_0$. It follows by definition that $\Phi_{\theta_0,\varepsilon}(\theta_0) = p_{\lambda,\varepsilon}(|\theta_0|)$. Furthermore, Lemma A.1 shows as in Proposition 3.1 that the even function $\Phi_{\theta_0,\varepsilon}(\theta) - p_{\lambda,\varepsilon}(|\theta|)$ is decreasing on $(0, |\theta_0|)$ and increasing on $(|\theta_0|, \infty)$, giving the desired result.

To prove part (b), it is sufficient to show that $|p_{\lambda,\varepsilon}(|\theta|) - p_\lambda(|\theta|)| \to 0$ uniformly on compact subsets of the parameter space as $\varepsilon \downarrow 0$. Since $p'_\lambda(\theta+)$ is nonincreasing on $(0,\infty)$,

$$|p_{\lambda,\varepsilon}(|\theta|) - p_\lambda(|\theta|)| \le \varepsilon \log\left[1 + \frac{|\theta|}{\varepsilon}\right] p'_\lambda(0+),$$

and because $p'_\lambda(0+) < \infty$, the right-hand side of the above inequality tends to 0 uniformly on compact subsets of the parameter space as $\varepsilon \downarrow 0$. □

Proof of Corollary 3.2. Let $\hat{\boldsymbol\beta}$ denote a maximizer of $Q(\boldsymbol\beta)$ and put $B = \{\boldsymbol\beta \in R^d : Q_{\varepsilon_0}(\boldsymbol\beta) \ge Q(\hat{\boldsymbol\beta})\}$ for some fixed $\varepsilon_0 > 0$. Then $B$ is compact and contains all $\hat{\boldsymbol\beta}_\varepsilon$ for $0 \le \varepsilon < \varepsilon_0$. Thus, Proposition 3.2 shows that for $\varepsilon < \varepsilon_0$,

$$\begin{aligned} Q(\hat{\boldsymbol\beta}) - Q(\hat{\boldsymbol\beta}_\varepsilon) &\le Q(\hat{\boldsymbol\beta}) - Q_\varepsilon(\hat{\boldsymbol\beta}) + Q_\varepsilon(\hat{\boldsymbol\beta}_\varepsilon) - Q(\hat{\boldsymbol\beta}_\varepsilon) \\ &\le |Q(\hat{\boldsymbol\beta}) - Q_\varepsilon(\hat{\boldsymbol\beta})| + |Q_\varepsilon(\hat{\boldsymbol\beta}_\varepsilon) - Q(\hat{\boldsymbol\beta}_\varepsilon)| \\ &\to 0. \end{aligned}$$

If $\boldsymbol\beta^*$ is a limit point of $\{\hat{\boldsymbol\beta}_\varepsilon\}_{\varepsilon\downarrow 0}$, then by the continuity of $Q(\boldsymbol\beta)$, $|Q(\boldsymbol\beta^*) - Q(\hat{\boldsymbol\beta})| = 0$ and so $\boldsymbol\beta^*$ is a maximizer of $Q(\boldsymbol\beta)$. □

Proof of Proposition 3.3. Given an initial value $\boldsymbol\beta^{(0)}$, let $\boldsymbol\beta^{(k)} = M^k(\boldsymbol\beta^{(0)})$ for $k \ge 1$; that is, $\{\boldsymbol\beta^{(k)}\}$ is the sequence of points that the MM algorithm generates starting from $\boldsymbol\beta^{(0)}$. Let $\Lambda$ denote the set of limit points of this sequence. For any $\boldsymbol\beta^* \in \Lambda$, passing to a subsequence we have $\boldsymbol\beta^{(k_n)} \to \boldsymbol\beta^*$. The quantity $Q_\varepsilon(\boldsymbol\beta^{(k_n)})$, since it is increasing in $n$ and bounded above, converges to a limit as $n \to \infty$. Thus, taking limits in the inequalities

$$Q_\varepsilon(\boldsymbol\beta^{(k_n)}) \le Q_\varepsilon\{M(\boldsymbol\beta^{(k_n)})\} \le Q_\varepsilon(\boldsymbol\beta^{(k_{n+1})})$$

gives $Q_\varepsilon(\boldsymbol{\beta}^*) = Q_\varepsilon\{\lim_{n\to\infty} M(\boldsymbol{\beta}_{k_n})\}$, assuming this limit exists. Of course, if $M(\boldsymbol{\beta})$ is continuous at $\boldsymbol{\beta}^*$, then we have $Q_\varepsilon(\boldsymbol{\beta}^*) = Q_\varepsilon\{M(\boldsymbol{\beta}^*)\}$, which implies that $\boldsymbol{\beta}^*$ is a stationary point of $Q_\varepsilon(\boldsymbol{\beta})$.

Note that $\Lambda$ is not necessarily nonempty in the above proof. However, we know that each $\boldsymbol{\beta}^{(k)}$ lies in the set $\{\boldsymbol{\beta} : Q_\varepsilon(\boldsymbol{\beta}) \geq Q_\varepsilon(\boldsymbol{\beta}_1)\}$, so if this set is compact, as is often the case, we may conclude that $\Lambda$ is indeed nonempty. □

**Acknowledgment.** The authors thank the referees for constructive suggestions.

Department of Statistics
Pennsylvania State University
University Park, Pennsylvania 16802-2111
USA
e-mail: dhunter@stat.psu.edu
e-mail: rli@stat.psu.edu